\documentclass[a4paper,10pt]{article}
\usepackage{amsfonts, amsthm}
\newcommand{\field}[1]{\mathbb{#1}}

\title{Surfaces obtained from $\field{C}P^{N-1}$ sigma models}

\author{\vspace{1cm}\\
         {\bf A. M.\ Grundland}$^{1,2}$
          \thanks{E-mail address:
           grundlan@crm.umontreal.ca}
             {\, and \,\bf \.{I}.\ Yurdu\c{s}en}$^1$
             \thanks{E-mail address:
       yurdusen@crm.umontreal.ca}
          \\
          \\$^1$Centre de Recherches Math\'{e}matiques, Universit\'{e} 
                   de Montr\'{e}al, 
               \\ CP 6128, Succ. Centre-Ville, Montr\'{e}al, 
                        Qu\'{e}bec H3C 3J7, Canada
                   \\
                 \\$^2$Universit\'{e} du Qu\'{e}bec, 
                     Trois-Rivi\`{e}res, CP500, QC, G9A 5H7, Canada}

\date{\today}
\begin{document}

\maketitle

\begin{abstract}
In this paper, the Weierstrass technique for harmonic maps
$S^2 \rightarrow \field{C}P^{N-1}$ is employed in order to obtain surfaces 
immersed in multidimensional Euclidean spaces. It is shown that 
if the $\field{C}P^{N-1}$ model equations are defined on the sphere 
$S^2$ and the associated action functional of this model is finite, then
the generalized Weierstrass formula for immersion describes conformally 
parametrized surfaces in the $su(N)$ algebra. In particular, 
for any holomorphic or antiholomorphic solution of this model 
the associated surface can be expressed in terms of 
an orthogonal projector of rank ($N-1$).
The implementation of this method is presented
for two-dimensional conformally parametrized surfaces immersed in 
the $su(3)$ algebra. 
The usefulness of the proposed approach is illustrated with 
examples, including the dilation-invariant meron-type solutions and 
the Veronese solutions for the $\field{C}P^{2}$ 
model. Depending on the location of the critical points (zeros and poles) of 
the first fundamental form associated with the meron solution, it is shown 
that the associated surfaces are semi-infinite cylinders. It is also 
demonstrated that surfaces related to holomorphic and mixed Veronese 
solutions are immersed in $\field{R}^{8}$ and $\field{R}^{3}$, 
respectively.
\end{abstract}




\section{Introduction \label{intro}}

The expression describing surfaces with zero mean curvature 
(i.e. minimal surfaces) which are immersed in three-dimensional 
Euclidean space was first formulated by A. Enneper \cite{Enneper} 
and K. Weierstrass \cite{Weierstrass} 
one and a half centuries ago. Since then this idea has been 
thoroughly generalized and developed (e.g. \cite{Bianchi, Dobriner, Darboux, Thomsen}). The subject was implemented 
by several authors (e.g. \cite{Fokas1, Friedrich, Budinich}) who produced 
several variants of the Weierstrass representation. For a 
comprehensive review of this topic see e.g. \cite{Bobenko1,Guest, Helein1, Helein2, Kenmotsu, Osserman} and 
references therein. 

More recently, this subject was 
substantially elaborated by B. Konopelchenko and I. Taimanov \cite{Konopelchenko}, 
who first established the Weierstrass formulae 
for any generic surface immersed in $\field{R}^{3}$. These formulae 
have been used extensively to study the global properties of surfaces 
in $\field{R}^{3}$, as well as their integrable deformations \cite{Konopelchenko2}. 
By simple analogy with surfaces in the 
$\field{R}^{3}$ case an extension of the Weierstrass procedure 
to multi-dimensional Euclidean and Riemannian spaces was proposed 
by B. Konopelchenko and G. Landolfi \cite{Konopelchenko3}. Their 
approach was successful for certain classes of conformally 
parametrized surfaces immersed in these spaces. However, this 
procedure has some limitation due to the assumption of a specific 
form of the Weierstrass system of $2N$ complex-valued functions 
which satisfy Dirac-type equations.

It was only in the past few years that the approach to the same 
problem was reformulated by exploiting the connection between 
generalized Weierstrass representations and the $\field{C}P^{1}$ 
sigma models, first established in $\field{R}^{3}$ \cite{Grundland1}. 
This idea allows one to generalize this connection 
for the $\field{C}P^{N-1}$ case and derive in the adjoint $SU(N)$ 
representation the corresponding moving frame of conformally 
parametrized surfaces in $\field{R}^{N^2-1}$ space. This modified 
Weierstrass representation \cite{Grundland2, Grundland3, Grundland4, Grundland5} has proven 
to be more general than the  one proposed in \cite{Grundland6} 
and to generate more diverse classes of surfaces (e.g. the Veronese 
surfaces). This algebraic description of surfaces on Lie groups 
and homogeneous spaces allows us to calculate some new expressions 
in closed form which determine the fundamental characteristics 
of these surfaces. For this purpose, using Cartan's 
language of moving frames we derive the structural equations 
for immersion (e.g. the fundamental forms, the Gaussian curvature and 
the mean curvature vector) for the $\field{C}P^{N-1}$ model. The 
$\field{C}P^{N-1}$ 
models have found many applications in physics, to such areas as 
two-dimensional gravity \cite{Gross}, string theory \cite{Polchinski}, 
quantum field theory \cite{Amit}, statistical physics \cite{Nelson} and 
fluid mechanics \cite{Chorin}. 

This paper is concerned with smooth, 
orientable two-dimensional surfaces immersed in 
multi-dimensional Euclidean spaces. The crux of the matter is that 
the equations determining the formula for immersion are formulated 
directly in terms of matrices which take their values in the Lie 
algebra $su(N)$. The main advantage of this procedure is that, in using an 
orthogonal projector satisfying the Euler-Lagrange equations of the 
given sigma model, it leads to simpler formulae and allows us to write the 
explicit form of some expressions which previously were too involved to be 
presented. 

The objective of this paper is to study certain geometrical aspects of 
surfaces associated with the $\field{C}P^{N-1}$ sigma models. 
In particular, we discuss in detail the necessary 
conditions for the existence of the radius vectors of surfaces 
associated with the $\field{C}P^{N-1}$ sigma model, 
which are expressed in terms of 
an orthogonal projector of rank $N-1$. Furthermore, we have shown that 
the Weierstrass formula for immersion of surfaces associated with mixed 
solutions of the $\field{C}P^{N-1}$ model is no longer proportional to a 
rank-one projector (unlike the case for holomorphic and 
antiholomorphic solutions). Next, it is demonstrated 
that a parametrized surface, related to a Veronese mixed solution (i.e. 
an extension of the holomorphic case) is immersed in three-dimensional 
Euclidean space. Finally, we construct meron-like solutions 
of the $\field{C}P^{2}$ model and determine their geometric characteristics. 

The plan of this paper is as follows. Section \ref{prilaminary} contains 
a brief account of basic definitions and properties concerning the 
$\field{C}P^{N-1}$ models and fixes the notation. We give a geometric 
formulation for the generalized Weierstrass formula for immersion of 
a surface $\mathcal{F}$ in $\field{R}^{N^2-1}$. Next, we show that 
if the $\field{C}P^{N-1}$ model is defined on the sphere 
$S^2$ and the corresponding action functional of this model is finite, then
a specific holomorphic function (corresponding to a component of the 
energy-momentum tensor of the $\field{C}P^{N-1}$ model) vanishes. 
In Section 3, we 
investigate in great detail the Veronese surfaces related to the 
$\field{C}P^{2}$ model and construct their geometric characteristics. We 
show that the holomorphic and mixed solutions are associated with surfaces 
immersed in $\field{R}^{8}$ and $\field{R}^{3}$, respectively. In Section 4, 
we discuss certain aspects of the projector formalism in the context of 
surfaces. In Section 5, we present examples of the application of our 
approach to the dilation-invariant solutions of meron type. 
We perform the analysis using quadratic differentials and calculate 
the geometric implications. Section 6 contains final remarks, 
identifies some open questions on the subject and 
proposes some possible future developments.


\section{Harmonic maps from $S^2$ to 
$\field{C}P^{N-1}$ and the Weierstrass representation
\label{prilaminary}}

This paper is devoted to the exploration of relations between the 
$\field{C}P^{N-1}$ sigma models and the generalized Weierstrass 
formula for the immersion of two-dimensional surfaces in 
multi-dimensional Euclidean 
spaces. To this end we briefly review some basic notions and properties 
of the $\field{C}P^{N-1}$ sigma models. For further details on this 
subject we refer the reader to e.g. \cite{Bobenko1, Guest, Helein1, Helein2, Zakrzewski, Zakharov} and references therein.

In studying the $\field{C}P^{N-1}$ models one is interested in maps 
of the form 
$[z]:\Omega \rightarrow \field{C}P^{N-1}$ 
(where $\Omega$ is an open, connected subset 
of a complex plane $\field{C}$) which are 
stationary points of the action functional \cite{Zakrzewski}
\begin{eqnarray}
S=\frac{1}{4}\int_{\Omega}(D_{\mu} z)^{\dagger} (D_{\mu} z) 
d\xi d\bar{\xi}\,, \qquad z^{\dagger}\cdot z=1\,, \nonumber \\
\field{C} \ni \xi=\xi^1+i \xi^2 \rightarrow 
z=(z_0, z_1, \ldots, z_{N-1}) \in \field{C}^N,
\label{action}
\end{eqnarray}
and thus are determined as solutions of the corresponding 
Euler-Lagrange equations. Here, 
$D_{\mu}$ denote covariant 
derivatives acting on $z:\Omega \rightarrow \field{C}^N$, 
defined by 
\begin{equation}
D_{\mu}z =\partial_{\mu}z- (z^{\dagger}\cdot \partial_{\mu}z)z
\,\, \in T_z S^{2N-1}\,,
\qquad \partial_{\mu}=\partial_{\xi^{\mu}}\,,
\quad \mu=1,2\,,
\label{covader}
\end{equation}
where $\xi$ and $\bar{\xi}$ are local coordinates in $\Omega$ and 
the symbol $\dagger$ denotes Hermitian conjugation. The covariant 
derivatives $D_{\mu}$ are orthogonal to the inhomogeneous coordinates 
$z$, since $z^{\dagger} D_{\mu} z =0 $ holds. They can be expressed 
in terms of a composite gauge field 
\begin{equation}
A_{\mu} = z^{\dagger} \partial_{\mu} z\,, \qquad A_{\mu}^{\dagger} = -A_{\mu}\,.
\label{compositegauge} 
\end{equation}
Here, $A_{\mu}$ is a pure imaginary function of $\xi^1$ and $\xi^2$. 
The action functional (\ref{action}) is invariant under global 
$U(N)$ transformations and also under the local $U(1)$ 
gauge transformation $z\, \rightarrow z^{\prime} = z e^{i\phi}$, 
where $\phi$ is a real-valued function. Note that the covariant 
derivatives $D_{\mu} z$ transform under the gauge transformation 
$D_{\mu} z\, \rightarrow D_{\mu} z^{\prime} = (D_{\mu} z) e^{i\phi}$, 
so that the dependence on the phase $\phi$ drops out of the action 
functional (\ref{action}) and so the model is really based on 
$\field{C}P^{N-1}$. In the homogeneous coordinates  
\begin{equation}
z={f}(f^{\dagger}\cdot f)^{-\frac{1}{2}}
\end{equation}
the equations of motion can be written in the form of a conservation law
\begin{equation}
\partial K - \bar{\partial} K^{\dagger}=0\,, 
\qquad
-i \partial K \in  su(N)\,,
\label{conservation}
\end{equation}
where $K$ and $K^{\dagger}$ are $N \times N$ matrices of the form 
\begin{eqnarray}
&&K=\frac{1}{f^{\dagger}\cdot f}
\left(\bar{\partial}f\otimes f^{\dagger}-f\otimes 
\bar{\partial}f^{\dagger}\right)
+\frac{f \otimes f^{\dagger}}{(f^{\dagger}\cdot f)^2}
\big(\bar{\partial}f^{\dagger}\cdot f - f^{\dagger}\cdot\bar{\partial}f
\big)\,, \nonumber \\
&&K^{\dagger}=\frac{1}{f^{\dagger}\cdot f}
\left(f\otimes \partial f^{\dagger}-\partial f\otimes 
f^{\dagger}\right)
+\frac{f \otimes f^{\dagger}}{(f^{\dagger}\cdot f)^2}
\big({\partial}f^{\dagger}\cdot f - f^{\dagger}\cdot {\partial}f
\big)\,.
\label{defmatK}
\end{eqnarray}
The symbols $\partial$ and $\bar{\partial}$ denote the standard derivatives 
with respect to $\xi$ and $\bar{\xi}$ respectively, {i.e.}
\begin{equation}
\partial=\frac{1}{2}\left(\partial_{\xi^1}-i\partial_{\xi^2}\right)\,, 
\qquad
\bar{\partial}=\frac{1}{2}\left(\partial_{\xi^1}+i\partial_{\xi^2}\right)\,.
\label{partials}
\end{equation}

Since the action (\ref{action}) is invariant under a global 
${U}(N)$ transformation, without loss of generality we can set one 
of the components of the vector field $f$ equal to $1$. Thus, in terms of 
these variables $f=(1, \bar{w}_1, \ldots, \bar{w}_N)^T$ the equations of 
motion for the $\field{C}P^{N-1}$ sigma model take the following form
\begin{eqnarray}
\partial\bar{\partial}w_i-\frac{2\bar{w}_i}{A_{N-1}}\partial w_i 
\bar{\partial} w_i-
\frac{1}{A_{N-1}}\sum_{j\ne i}^{N-1} \bar{w}_j (\partial w_i 
\bar{\partial} w_j + \bar{\partial} w_i \partial w_j)=0\,,\nonumber\\
\partial\bar{\partial}\bar{w}_i-\frac{2 w_i}{A_{N-1}}\partial 
\bar{w}_i \bar{\partial} \bar{w}_i-
\frac{1}{A_{N-1}} \sum_{j\ne i}^{N-1} w_j (\partial 
\bar{w}_i \bar{\partial} \bar{w}_j + \bar{\partial} \bar{w}_i 
\partial \bar{w}_j)=0\,, \label{cpn}
\end{eqnarray}
where $i=1,2,...,N-1$ and $A_{N-1}=1+\sum_i^{N-1}w_i\bar{w}_i$. 
In what follows we refer to (\ref{cpn}) as the equations of the $\field{C}P^{N-1}$ sigma model.

It is instructive to express the Euler-Lagrange equations using the 
$N \times N$ orthogonal projector $P$ of rank ($N-1$) defined 
on the orthogonal complement to the complex line in 
$\field{C}^{N}$, 
\begin{equation}
P=I_{N}-\frac{f \otimes f^{\dagger}}{f^{\dagger}\cdot f}\,,
\qquad
P^{\dagger}=P\,, \qquad P^2=P\,,
\label{projector}
\end{equation}
where $I_N$ is the $N \times N$ identity matrix.
Hence, the Euler-Lagrange equation (\ref{conservation}) takes the simpler 
form
\begin{equation}
\partial [\bar{\partial} P, P] + \bar{\partial} [\partial P, P]=0\,.
\label{conservation2}
\end{equation}

After expressing the Euler-Lagrange equations (\ref{conservation2}) 
as a conservation law, we are able to formulate the Weierstrass formula 
for the immersion of two-dimensional surfaces in multi-dimensional Euclidean 
space. Based on Poincar\'e's lemma, there exists a closed matrix-valued 
$1$-form, 
\begin{equation}
dX=i(-[\partial P, P]d\xi + [\bar{\partial}P, P] d\bar{\xi})\,. 
\label{complex1formimmer}
\end{equation}
From the closure of the $1$-form $dX$ (i.e. $d(dX)=0$) it follows that 
the integral
\begin{equation}
X(\xi, \bar{\xi})=
i\int_{\gamma} (-[\partial P, P]d\xi + [\bar{\partial}P, P] d\bar{\xi})\,,
\label{intforimm}
\end{equation}
depends only on the end points of the curve $\gamma$ (i.e. it is locally 
independent of the trajectory in $\field{C}$). Note that 
(\ref{conservation2}) is invariant under the conformal transformation 
(i.e. the change of independent variables $\xi \rightarrow \alpha (\xi)$ 
and $\bar{\xi} \rightarrow \bar{\alpha} (\bar{\xi})$). Such a transformation 
establishes a reparametrization of the surface $\mathcal{F}$ written 
in terms of an integral of a $1$-form (\ref{intforimm}) which remains 
the same geometrical object.

For the analytical description of a two-dimensional surface $\mathcal{F}$ 
it is convenient to use the Lie algebra isomorphism and identify 
the ($N^2-1$)-dimensional Euclidean space with the $su(N)$ algebra 
\begin{equation}
\field{R}^{N^2-1} \simeq {su}(N)\,.  
\label{isomorphism}
\end{equation}
For uniformity we use the scalar product on ${su}(N)$ in the form 
\begin{equation}
<A, B>=-\frac{1}{2} {\rm tr}(AB)\,,
\qquad
A,B \in su(N)\,,
\end{equation}
rather than the Killing form of $su(N)$ given by the formula
\begin{eqnarray}
\mathcal{B} (A,B) = 2\, N {\rm tr} (AB)\,,
\end{eqnarray}
which is negative definite \cite{Helgason}. Consequently the 
first fundamental form $I$ is given by \cite{ghy}
\begin{equation}
I=-Jd\xi^2 +\frac{2}{f^{\dagger}\cdot f} \bar{\partial}f^{\dagger} P \partial f
d\xi d\bar{\xi} - \bar{J} d{\bar{\xi}}^2\,, 
\label{genfirstfundform}
\end{equation}
where the complex-valued functions $J$ and $\bar{J}$ are given by 
\begin{equation}
J=\frac{1}{f^{\dagger}\cdot f}\partial f^{\dagger} P \partial f\,,
\qquad
\bar{J}=\frac{1}{f^{\dagger}\cdot f}\bar{\partial} f^{\dagger} 
P \bar{\partial} f\,. \label{defJ}
\end{equation}
They satisfy
\begin{equation}
\bar{\partial} J=0\,, \qquad \partial \bar{J}=0\,,
\end{equation}
whenever $f$ is a solution of the equations of motion (\ref{conservation}). 
The quantities $J$ and $\bar{J}$ are invariant under the global ${U}(N)$ 
transformation, {i.e.,} $f \rightarrow a f$, $a \in {U}(N)$. 
From the physical point of view, $J=({D}z)^{\dagger}\cdot Dz$ is related to the 
energy-momentum tensor \cite{Zakrzewski}.

The integral representation (\ref{intforimm}) defines a mapping 
$X: \Omega \ni (\xi, \bar{\xi}) \rightarrow X (\xi, \bar{\xi})
\in su(N)$. We treat each element of the real-valued $su(N)$ matrix function 
$X$ as coordinates of a two-dimensional surface $\mathcal{F}$ immersed in 
$\field{R}^{N^2-1}$. This map $X$ is called the generalized Weierstrass formula 
for immersion. The projector $P$ is invariant under the transformation 
$P \rightarrow UPU^{\dagger}$, where $U\in U(N)$ and thus the geometry 
of the surface $\mathcal{F}$ associated with a solution of 
(\ref{conservation2}) admits the symmetry equivalence class of solutions 
of (\ref{conservation2}). In this setting, our generalization lies in the realization 
that most of the properties of the associated surfaces with the 
$\field{C}P^{N-1}$ sigma models can be described using an orthogonal 
projector. The complex tangent vectors of this immersion 
are 
\begin{equation}
\partial X = i K^{\dagger}\,,
\qquad
\bar{\partial} X = i K,
\end{equation}
where we use (\ref{defmatK})
\begin{equation}
K=[\bar{\partial} P, P]\,,
\qquad
K^{\dagger} = -[\partial P, P]\,.
\label{defmatKincommutator}
\end{equation}

From the conservation law (\ref{conservation}), it is convenient 
to decompose the matrix $K$ as follows
\begin{equation}
K=M+L\,,
\label{decomofK}
\end{equation}
where 
\begin{equation}
M=(I_{N}-P) \bar{\partial} P\,,
\qquad
L=-\bar{\partial} P (I_{N}-P)\,.
\label{decomofKplusML}
\end{equation}
It was shown in \cite{Grundland6} that the matrices $M$ and $L$ satisfy 
the same conservation law (\ref{conservation}) as the matrix 
$K$
\begin{equation}
\partial M = \bar{\partial} M^{\dagger}\,,
\qquad
\partial L = \bar{\partial} L^{\dagger}\,,
\label{sameconser}
\end{equation}
and the matrices $M$ and $L$ differ by a total divergence
\begin{equation}
M=L+\bar{\partial} P\,.
\label{differtotal}
\end{equation}

Let us now discuss the existence of certain classes of surfaces immersed 
in the $su(N)$ algebra under the hypotheses that the $\field{C}P^{N-1}$ 
model is defined on the sphere $S^2$ and its corresponding action 
functional (\ref{action}) is finite. In this 
case, the procedure for constructing the general class of solutions 
of the Euclidean two-dimensional $\field{C}P^{N-1}$ model 
(\ref{cpn}) was derived by A. Din and W. Zakrzewski \cite{Din} 
and R. Sasaki \cite{Sasaki}. As a 
result, one gets three classes of solutions, namely (i) holomorphic 
(i.e. $\bar{\partial} f = 0$), (ii) antiholomorphic (i.e. ${\partial} f = 0$) 
and (iii) mixed. The mixed solutions can be determined from either 
the holomorphic or the antiholomorphic nonconstant functions by the following 
procedure. The successive application, say $k$ times with $k \le N-1$, of the 
operator $P_{+}$ defined by its action on vector-valued functions on 
$\field{C}^{N}$ \cite{Zakrzewski}
\begin{eqnarray}
P_{+}: f \in \field{C}^N \rightarrow
P_{+}f=\partial f -f \frac{f^{\dagger} \partial f}{f^{\dagger} f}\,,
\quad \bar{\partial} f = 0\,,
\label{operatorintermsoff}
\end{eqnarray}
starting from any nonconstant holomorphic function $f\in \field{C}^N$, 
allows one to find mixed solutions 
\begin{equation}
f^k=P_{+}^k f\,, \qquad k=0,1, \ldots, N-1\,,
\end{equation}
which represent harmonic maps from $S^2$ to the $\field{C}P^{N-1}$ 
sigma model. Here, $P_{+}^0 = id$.

Note that the holomorphic function $f\in \field{C}^N$, used 
in (\ref{operatorintermsoff}), could be replaced by any nonconstant 
antiholomorphic function. The mixed solutions $f^k$ are constructed 
in the same way, except that the derivative $\partial$ is replaced by 
$\bar{\partial}$ in the definition of the operator $P_{-}$. Thus, we have 
\begin{eqnarray}
P_{-}f=\bar{\partial} f -f \frac{f^{\dagger} \bar{\partial} f}
{f^{\dagger} f}\,,
\quad \partial f = 0\,,
\label{operatorintermsoffforantihol}
\end{eqnarray}
which yields complementary results.

Under the above hypotheses, the considered surfaces are conformally 
parametrized and the 
first fundamental form (\ref{genfirstfundform}) becomes
\begin{equation}
I=\frac{2}{f^{\dagger}\cdot f} \bar{\partial}f^{\dagger} P \partial f
d\xi d\bar{\xi}\,. 
\label{conformalfirstfundform}
\end{equation}
In order to demonstrate that the complex-valued functions 
$J$ and $\bar{J}$ vanish it is sufficient to consider
the orthogonality relation 
\begin{equation}
(P_{+}^i f)^{\dagger} \cdot P_{+}^j f = 0\,, 
\qquad
i \ne j\,.
\label{orthogonalitypr}
\end{equation}
for $i=k$ and $j=k+2$ with arbitrary $k=0,1, \ldots, N-1$. 
Denoting  $\widetilde{f} = P_{+}^k f $, we get
\begin{eqnarray}
0= \widetilde{f}^{\dagger} \cdot (P_{+}^2 \widetilde{f}) &=& 
\widetilde{f}^{\dagger} \cdot \left( 
\partial (P_{+} \widetilde{f}) -(P_{+} \widetilde{f}) \frac{(P_{+} 
\widetilde{f})^{\dagger} \partial (P_{+} \widetilde{f})}{(P_{+} 
\widetilde{f})^{\dagger} (P_{+} \widetilde{f})}
\right) \nonumber \\
&=& \widetilde{f}^{\dagger} \cdot \partial (P_{+} \widetilde{f})
\nonumber \\
&=& -\partial \widetilde{f}^{\dagger} \cdot (P_{+} \widetilde{f})\,,
\label{proofderivation}
\end{eqnarray}
where for the last two equalities we used the orthogonality condition 
$\widetilde{f}^{\dagger} \cdot P_{+} \widetilde{f} = 0$. The right 
hand side of the last equality in (\ref{proofderivation}) can also be 
written in terms of the complex-valued functions $J$ and $\bar{J}$ given in 
(\ref{defJ})
\begin{eqnarray}
0= - \partial \widetilde{f}^{\dagger} \cdot \left(
\partial \widetilde{f} - \frac{\widetilde{f} 
\otimes \widetilde{f}^{\dagger}}{\widetilde{f}^{\dagger} 
\cdot \widetilde{f}} \partial \widetilde{f}
\right)
=-(\widetilde{f}^{\dagger} \cdot \widetilde{f}) \widetilde{J}\,.
\label{showingj0}
\end{eqnarray}
Since $\widetilde{f}^{\dagger} \cdot \widetilde{f} \ne 0$, we have 
$\widetilde{J}=0$. Note that for the holomorphic and antiholomorphic 
solutions $f$ of the $\field{C}P^{N-1}$ model equations 
(\ref{conservation}) the corresponding complex-valued functions 
$J$ and $\bar{J}$, given in (\ref{defJ}), vanish identically. 

In particular one can present an analogue of the Bonnet 
theorem. If we consider the holomorphic or antiholomorphic solutions 
of the $\field{C}P^{N-1}$ model under the above hypotheses, the Weierstrass 
formula for immersion $X$ of a surface $\mathcal{F}$ 
can be expressed in terms of 
the orthogonal projector of rank $(N-1)$ by the formula
\begin{equation}
X(\xi, \bar{\xi})= \epsilon i \left(\frac{1-N}{N} I_{N} + P 
\right)\,,
\qquad
\epsilon=\pm 1\,.
\label{radiuswithepsilonrevised}
\end{equation}
The surface $\mathcal{F}$ is determined uniquely up to Euclidean motions 
by its first and second fundamental forms
\begin{equation}
I={\rm tr} (\partial P \bar{\partial} P) d\xi d\bar{\xi}\,,
\label{firstfunwithp}
\end{equation}
and
\begin{eqnarray}
\!\!\!\!II\!\!=\!\!\epsilon i \Big\{\!\!
(\partial^2 P - \Gamma_{11}^1 \partial P - \Gamma_{11}^2 \bar{\partial} P) d\xi^2 \!+\! 
2 \partial \bar{\partial} P d\xi d\bar{\xi} \!+\!
(\bar{\partial}^2 P - \Gamma_{22}^1 \partial P - \Gamma_{22}^2 \bar{\partial} P) d\bar{\xi}^2
\!\!\Big\}\!,
\label{secondfunwithp}
\end{eqnarray}
respectively, where the Christoffel symbols of the second kind are given by 
\begin{eqnarray}
\Gamma_{11}^1 = \frac{{\rm tr} (\partial^2 P \bar{\partial} P)}{{\rm tr} (\partial P \bar{\partial} P)}\,,
\quad
\Gamma_{11}^2 = \frac{{\rm tr} (\partial^2 P \partial P)}{{\rm tr} (\partial P \bar{\partial} P)}\,,
\quad
\Gamma_{12}^1 = \Gamma_{21}^1 = 0\,, \nonumber \\
\Gamma_{22}^1 = \frac{{\rm tr} (\bar{\partial}^2 P \bar{\partial} P)}{{\rm tr} (\partial P \bar{\partial} P)}\,,
\quad
\Gamma_{22}^2 = \frac{{\rm tr} (\bar{\partial}^2 P \partial P)}{{\rm tr} (\partial P \bar{\partial} P)}\,,
\quad
\Gamma_{12}^2 = \Gamma_{21}^2 = 0\,.
\label{Christoffelsymbolswithp}
\end{eqnarray}
In the case of holomorphic or antiholomorphic solutions $f$ of the 
$\field{C}P^{N-1}$ model, according to \cite{ghy}, the matrix $K$ 
can be expressed in the simple form
\begin{equation}
K = \epsilon \bar{\partial} P\,,
\qquad
K^{\dagger} = \epsilon \partial P\,,
\qquad
\epsilon = \pm 1\,.
\label{defKfornewproposition}
\end{equation}



\section{Veronese surfaces for the $\field{C}P^{2}$ model}

One of the simplest applications of a result concerning solutions 
of the $\field{C}P^{N-1}$ sigma model (\ref{cpn}) is the 
Veronese sequence \cite{Bolton}
\begin{eqnarray}
f=\left(1, \sqrt{\left(\begin{array}{c}
N-1  \\
1  \\
\end{array}\right)} \,\,\xi , 
\ldots, \sqrt{\left(\begin{array}{c}
N-1  \\
r  \\
\end{array}\right)} \,\,\xi^r ,
\ldots,  \xi^{N-1}
\right)\,.
\label{genveronese}
\end{eqnarray}
For all of the above Veronese solutions the first fundamental form 
is conformal and given by 
\begin{equation}
I=(N-1) (1+|\xi|^2)^{-2} d \xi d \bar{\xi}\,. 
\label{firstfundamentalforforgeneralveronese}
\end{equation}
Since $g_{11}=g_{22}=0$, the Gaussian curvature is computed from 
the following formula \cite{ghy}
\begin{equation}
\mathcal{K}=-(g_{12})^{-1}\bar{\partial}\partial\ln g_{12}\,,
\label{Gausscurvatureforforgeneralveronese1}
\end{equation}
and for the Veronese solutions it is found to be
\begin{equation}
\mathcal{K}=\frac{4}{N-1}\,.
\label{Gausscurvatureforforgeneralveronese2}
\end{equation}

From now on we will only be concerned with the 
$\field{C}P^{2}$ model ($N=3$) for which (\ref{cpn}) 
becomes
\begin{eqnarray}
&&\partial\bar{\partial}w_1-\frac{2\bar{w}_1}{A_2}\partial w_1 
\bar{\partial} w_1-\frac{\bar{w}_2}{A_2}(\partial w_1 \bar{\partial} 
w_2 + \bar{\partial} w_1 \partial w_2)=0\,,\nonumber\\
&&\partial\bar{\partial}w_2-\frac{2\bar{w}_2}{A_2}\partial w_2 
\bar{\partial} w_2-\frac{\bar{w}_1}{A_2}(\partial w_1 \bar{\partial} 
w_2 + \bar{\partial} w_1 \partial w_2)=0\,, \nonumber\\
&&A_2=1+w_1\bar{w}_1+w_2\bar{w}_2\,,
\label{cp2}
\end{eqnarray}
together with their complex conjugate equations. The Veronese vector $f$ for this model 
is given by 
\begin{equation}
f=(1, \sqrt{2}\,\xi, \xi^2)\,.
\label{cp2veronese}
\end{equation}

The method for finding the radius vector $\vec{X}$ through 
the use of the generalized Weierstrass formula for immersion of $2D$ 
surfaces in $\field{R}^8$ was proposed in \cite{Grundland2} 
and \cite{ghy}. According to \cite{ghy}, the real components 
of the corresponding $1$-forms for any solution of the 
$\field{C}P^{2}$ model are 
\begin{eqnarray}
dX_1&=&\frac{1}{2{A_2}^{\!\!2}}\Big(\big[(w_2^2-w_1^2)(\bar{w}_1\partial 
\bar{w}_2-\bar{w}_2\partial \bar{w}_1)-(\bar{w}_2^2-\bar{w}_1^2)({w}_1
\partial {w}_2-{w}_2\partial {w}_1) \nonumber \\
&\,&-w_2\partial \bar{w}_1+\bar{w}_2\partial {w}_1-w_1\partial \bar{w}_2
+\bar{w}_1\partial {w}_2\big] d\xi + \rm{c.c.} \Big),
\nonumber \\
dX_2&=&\frac{i}{2{A_2}^{\!\!2}}\Big(\big[(w_1^2+w_2^2)(\bar{w}_2\partial 
\bar{w}_1-\bar{w}_1\partial \bar{w}_2)+(\bar{w}_1^2+\bar{w}_2^2)({w}_2
\partial {w}_1-{w}_1\partial {w}_2) \nonumber \\
&\,&+w_2\partial \bar{w}_1+\bar{w}_2\partial {w}_1-w_1\partial \bar{w}_2-
\bar{w}_1\partial {w}_2\big] d\xi - \rm{c.c.} \Big),
\nonumber \\
dX_3&=&\frac{1}{2{A_2}^{\!\!2}}\Big(\big[w_2\partial \bar{w}_2-{w}_1\partial 
\bar{w}_1-\bar{w}_2\partial {w}_2+\bar{w}_1\partial {w}_1 
\nonumber \\
&\,&+ 2 |w_1|^2 ({w}_2\partial \bar{w}_2-\bar{w}_2\partial {w}_2) 
- 2 |w_2|^2 ({w}_1\partial \bar{w}_1-\bar{w}_1\partial {w}_1) \big] d\xi + 
\rm{c.c.} \Big),
\nonumber \\
dX_4&=&\frac{\sqrt{3}}{2{A_2}^{\!\!2}}\Big(\big[w_1\partial \bar{w}_1+{w}_2
\partial \bar{w}_2-\bar{w}_1\partial {w}_1-\bar{w}_2\partial {w}_2 \big] 
d\xi + \rm{c.c.} \Big),
\nonumber \\
dX_5&=&-\frac{i}{2{A_2}^{\!\!2}}\Big(\big[(1+\bar{w}_1^2+|w_2|^2)
\partial w_1 + (1+{w}_1^2+|w_2|^2) \partial \bar{w}_1
\nonumber \\
&\,&+({w}_2\partial \bar{w}_2-\bar{w}_2\partial {w}_2)(w_1-\bar{w}_1)\big] 
d\xi - \rm{c.c.} \Big),
\nonumber \\
dX_6&=&-\frac{i}{2{A_2}^{\!\!2}}\Big(\big[(1+\bar{w}_2^2+|w_1|^2)
\partial w_2 + (1+{w}_2^2+|w_1|^2) \partial \bar{w}_2
\nonumber \\
&\,&+({w}_1\partial \bar{w}_1-\bar{w}_1\partial {w}_1)(w_2-\bar{w}_2)\big] 
d\xi - \rm{c.c.} \Big),
\nonumber \\
dX_7&=&\frac{1}{2{A_2}^{\!\!2}}\Big(\big[(1-{w}_1^2+|w_2|^2)\partial 
\bar{w}_1 - (1-\bar{w}_1^2+|w_2|^2) \partial {w}_1
\nonumber \\
&\,&+(\bar{w}_2\partial {w}_2-{w}_2\partial \bar{w}_2)(w_1+\bar{w}_1)\big] 
d\xi + \rm{c.c.} \Big),
\nonumber \\
dX_8&=&\frac{1}{2{A_2}^{\!\!2}}\Big(\big[(1-{w}_2^2+|w_1|^2)\partial 
\bar{w}_2 - (1-\bar{w}_2^2+|w_1|^2) \partial {w}_2
\nonumber \\
&\,&+(\bar{w}_1\partial {w}_1-{w}_1\partial \bar{w}_1)(w_2+\bar{w}_2)\big] 
d\xi + \rm{c.c.} \Big)\,. 
\label{weiestrassforcp2}
\end{eqnarray}
For any holomorphic solution ($w_1, w_2$) of the $\field{C}P^{2}$ model 
the above $8$ real-valued $1$-forms can easily be integrated to give
the components of the radius vector
\begin{equation}
\vec{X}(\xi, \bar{\xi})=\Big(X_1(\xi, \bar{\xi}),\dots, X_8(\xi, \bar{\xi})\Big)\,,
\label{radvec}
\end{equation}
of a two-dimensional surface in $\field{R}^8$
\begin{eqnarray}
&&X_1=\frac{w_1\bar{w}_2+\bar{w}_1w_2}{2\,A_2}\,, \qquad 
X_2=i\frac{w_1\bar{w}_2-\bar{w}_1w_2}{2\,A_2}\,, \qquad 
X_3=\frac{|w_1|^2-|w_2|^2}{2\,A_2}\,, \nonumber \\
&&X_4=-\sqrt{3}\frac{|w_1|^2+|w_2|^2}{2\,A_2}\,, \qquad 
X_5=-i\frac{w_1-\bar{w}_1}{2\,A_2}\,, \qquad
X_6=-i\frac{w_2-\bar{w}_2}{2\,A_2}\,, \nonumber \\
&&X_7=-\frac{w_1+\bar{w}_1}{2\,A_2}\,, \qquad
X_8=-\frac{w_2+\bar{w}_2}{2\,A_2}\,,
\label{coordinatesofrcp2}
\end{eqnarray}
where we choose the integration constants to be zero.

Hence, using the Weierstrass formula for 
immersion (\ref{coordinatesofrcp2}) we obtain that the radius vector $\vec{X}$ 
of a two-dimensional parametrized surface (for the Veronese solution (\ref{cp2veronese})) 
is immersed in $\field{R}^8$. Its components are 
\begin{eqnarray}
&&X_1=\frac{|\xi|^2 (\xi+\bar{\xi})}{\sqrt{2}\,(1+|\xi|^2)^2}=\frac{\sqrt{2}\,x(x^2+y^2)}{(1+x^2+y^2)^2}\,, \nonumber \\
&&X_2=-i\frac{|\xi|^2 (\xi-\bar{\xi})}{\sqrt{2}\,(1+|\xi|^2)^2}=\frac{\sqrt{2}\,y(x^2+y^2)}{(1+x^2+y^2)^2}\,, \nonumber \\
&&X_3=-\frac{|\xi|^2(|\xi|^2-2)}{2 (1+|\xi|^2)^2}=-\frac{(x^2+y^2)(x^2+y^2-2)}{2 (1+x^2+y^2)^2}\,, \nonumber \\
&&X_4=-\sqrt{3}\frac{|\xi|^2(|\xi|^2+2)}{2 (1+|\xi|^2)^2}=-\sqrt{3}\frac{(x^2+y^2)(x^2+y^2+2)}{2 (1+x^2+y^2)^2}\,,\nonumber \\
&&X_5=-i\frac{\xi-\bar{\xi}}{\sqrt{2}\,(1+|\xi|^2)^2}=\frac{\sqrt{2}\,y}{(1+x^2+y^2)^2}\,, \nonumber \\
&&X_6=-i\frac{\xi^2-\bar{\xi}^2}{{2}\,(1+|\xi|^2)^2}=\frac{{2}\,xy}{(1+x^2+y^2)^2}\,,\nonumber \\
&&X_7=-\frac{\xi+\bar{\xi}}{\sqrt{2}\,(1+|\xi|^2)^2}=-\frac{\sqrt{2}\,x}{(1+x^2+y^2)^2}\,,\nonumber \\
&&X_8=-\frac{\xi^2+\bar{\xi}^2}{{2}\,(1+|\xi|^2)^2}=\frac{-x^2+y^2}{(1+x^2+y^2)^2}\,, 
\label{compr81}
\end{eqnarray}
where we used $\xi=x+iy$. The components $X_i$ ($i=1, \ldots, 8$) given in 
(\ref{compr81}) satisfy the equation of an affine sphere 
\begin{equation}
4X_1^2+4X_2^2+4X_3^2+\frac{2}{\sqrt{3}}X_4+X_5^2+X_6^2+X_7^2+X_8^2=0\,.
\label{surfaceeqforcp2holomorphic}
\end{equation}

We can now proceed to construct a mixed solution which, as is well-known \cite{Zakrzewski}, 
can be obtained directly from the holomorphic one. 
Applying the operator $P_{+}$, given by (\ref{operatorintermsoff}), to 
the vector field (\ref{cp2veronese}), we obtain the mixed solution in 
the form
\begin{equation}
P_{+}f=\frac{\sqrt{2}}{1+|\xi|^2}(-\sqrt{2}\,\bar{\xi}, 1-|\xi|^2, \sqrt{2}\,{\xi})\,.
\label{newsol}
\end{equation}
Let us note that for the $\field{C}P^{2}$ model, the repeated applications 
of the operator $P_{+}$ to a holomorphic solution $f$ only lead to the mixed 
solution (\ref{newsol}) and an antiholomorphic one $P_{+}^2 f$, 
since $P_{+}^3 f = 0$. 
Thus, the holomorphic and mixed solutions considered here indeed 
constitute a complete set of solutions for the $\field{C}P^{2}$ model. 
Using $U(1)$ invariance of the $\field{C}P^{2}$ model we can normalize 
(\ref{newsol}) to the following vector 
\begin{equation}
f_1=(1, \widetilde{w_1}, \widetilde{w_2})\,,
\label{cp2veronesemix}
\end{equation}
where we denote
\begin{equation}
\widetilde{w_1}=\frac{|\xi|^2-1}{\sqrt{2}\, \bar{\xi}}\,,
\qquad
\widetilde{w_2}=-\frac{\xi}{\bar{\xi}}\,.
\label{mixedsolcp2}
\end{equation} 
Then, substituting (\ref{mixedsolcp2}) into (\ref{weiestrassforcp2}) 
and integrating, we obtain a two-dimensional parametrized 
surface immersed in $\field{R}^3$
\begin{eqnarray}
&&X_1=-X_7=\frac{\xi+\bar{\xi}}{\sqrt{2}\,(1+|\xi|^2)}=\frac{\sqrt{2}\,x}{1+x^2+y^2}\,,  \nonumber \\
&&X_3=\frac{X_4}{\sqrt{3}}=\frac{1}{1+|\xi|^2}=\frac{1}{1+x^2+y^2}\,, \nonumber \\
&&X_2=X_5=-i\frac{\xi-\bar{\xi}}{\sqrt{2}\,(1+|\xi|^2)}=\frac{\sqrt{2}\,y}{1+x^2+y^2}\,,  \nonumber \\
&&X_6=X_8=0\,. 
\label{radiusvecmix}
\end{eqnarray}
Note that the components of the radius vector $\vec{X}$ 
in (\ref{radiusvecmix}) satisfy the following relation
\begin{equation}
X_1^2+X_2^2+(\sqrt{2}X_3-\frac{1}{\sqrt{2}})^2=\frac{1}{2}\,.
\label{eqsurmix}
\end{equation}
Equation (\ref{eqsurmix}) represents an ellipsoid, centered at the point $(0,0,\frac{1}{2})$ 
in $\field{R}^3$. So, this case corresponds to the immersion of the $\field{C}P^{2}$ model 
into the $\field{C}P^{1}$ model.

Let us now explore some geometrical characteristics of surfaces corresponding to two different 
solutions of the $\field{C}P^{2}$ model. In the holomorphic case (\ref{cp2veronese}) the 
orthogonal projector has the following form 
\begin{equation}
P=\frac{1}{(1+|\xi|^2)^2} \left(\begin{array}{ccc}
|\xi|^2(2+|\xi|^2) & -\sqrt{2}\,\xi & -\xi^2 \\
-\sqrt{2}\,\bar{\xi} & 1+|\xi|^4 & -\sqrt{2}\,|\xi|^2 \xi \\
-\bar{\xi}^2 & -\sqrt{2}\,|\xi|^2 \bar{\xi} & 1+2|\xi|^2
\end{array}
\right)\,, \label{projectorcp2}
\end{equation}
where rank$P$ $=$ $2$ and tr$P$ $=2$. 
The surface is determined by 
(\ref{compr81}) and its induced metric is conformal 
\begin{equation}
g_{11}=g_{22}=0, \qquad g_{12}=\frac{1}{(1+|\xi|^2)^2}\,.
\label{indmetricforholomorphicus}
\end{equation}
The nonzero Christoffel symbols of the second kind are
\begin{equation}
\Gamma_{11}^1=-\frac{2\bar{\xi}}{1+|\xi|^2}\,, \qquad
\Gamma_{22}^2=-\frac{2 \xi}{1+|\xi|^2}\,.
\label{Christoffelforholomorphicus}
\end{equation}
The first fundamental form and the Gaussian curvature are given by
\begin{equation}
I=\frac{2}{(1+|\xi|^2)^2}d \xi d \bar{\xi}\,, \qquad \mathcal{K}=2,
\label{firstGaussforholomorphicus}
\end{equation}
respectively. Making use of the expression (\ref{compr81}) for the radius 
vector $\vec{X}$ we can explicitly write the second fundamental form $II$ of the 
surface in the equivalent matrix form. The components of the matrix $II$ are
\begin{eqnarray}
&&II_{11}=\frac{2i}{(1+|\xi|^2)^4}(\bar{\xi}^2 d\xi^2 + (4|\xi|^2-2) d\xi d\bar{\xi} + \xi^2 d\bar{\xi}^2)\,, \nonumber \\
&&II_{12}=\frac{2 \sqrt{2} i}{(1+|\xi|^2)^4}(-\bar{\xi} d\xi^2 + 2\xi(|\xi|^2-2) d\xi d\bar{\xi} + \xi^3 d\bar{\xi}^2)\,, \nonumber \\
&&II_{13}=\frac{2i}{(1+|\xi|^2)^4}(d\xi^2 - 6\xi^2 d\xi d\bar{\xi} + \xi^4 d\bar{\xi}^2)\,, \nonumber \\
&&II_{21}=\frac{2 \sqrt{2} i}{(1+|\xi|^2)^4}(\bar{\xi}^3 d\xi^2 + 2\bar{\xi}(|\xi|^2-2) d\xi d\bar{\xi} - \xi d\bar{\xi}^2)\,, \nonumber \\
&&II_{22}=\frac{4 i}{(1+|\xi|^2)^4}(-\bar{\xi}^2 d\xi^2 + (1+|\xi|^4-4|\xi|^2) d\xi d\bar{\xi} - \xi^2 d\bar{\xi}^2)\,, \nonumber \\
&&II_{23}=\frac{2 \sqrt{2} i}{(1+|\xi|^2)^4}(\bar{\xi} d\xi^2 + 2\xi(1-2|\xi|^2) d\xi d\bar{\xi} - \xi^3 d\bar{\xi}^2)\,, \nonumber \\
&&II_{31}=\frac{2i}{(1+|\xi|^2)^4}(\bar{\xi}^4 d\xi^2 - 6\bar{\xi}^2 d\xi d\bar{\xi} + d\bar{\xi}^2)\,, \nonumber \\
&&II_{32}=\frac{2 \sqrt{2} i}{(1+|\xi|^2)^4}(-\bar{\xi}^3 d\xi^2 + 2\bar{\xi}(1-2|\xi|^2) d\xi d\bar{\xi} + \xi d\bar{\xi}^2)\,, \nonumber \\
&&II_{33}=\frac{2i}{(1+|\xi|^2)^4}(\bar{\xi}^2 d\xi^2 + 2|\xi|^2(2-|\xi|^2) d\xi d\bar{\xi} + \xi^2 d\bar{\xi}^2)\,.
\label{secfuncomp}
\end{eqnarray}
The mean curvature $\mathcal{H}={\partial \bar{\partial} X}/{g_{12}}$, 
written as a matrix, takes the form
\begin{equation}
\mathcal{H}=\frac{4i}{(1+|\xi|^2)^2} \left(\begin{array}{ccc}
2|\xi|^2-1 & \sqrt{2}\,\xi (|\xi|^2-2) & -3\xi^2 \\
\sqrt{2}\,\bar{\xi} (|\xi|^2-2) & 1+|\xi|^2(|\xi|^2-4) & -\sqrt{2}\, \xi (2|\xi|^2-1)  \\
-3\bar{\xi}^2 & -\sqrt{2}\,\bar{\xi} (2|\xi|^2-1)  & -|\xi|^2(|\xi|^2-2)
\end{array}
\right)\,, 
\label{meancurh}
\end{equation}
where rank $\mathcal{H}=2$ and tr $\mathcal{H}=0$.
The total energy \cite{Zakrzewski} for the holomorphic solution 
(\ref{cp2veronese}) is finite over all space 
\begin{equation}
u=\ln\left(\frac{|\partial w_1|^2+|\partial w_2|^2+|w_2\partial w_1-
w_1\partial w_2|^2}{{A_2}^{2}}\right)=\ln\left(\frac{2}{{(1+|\xi|^2)}^{2}}\right)\,.
\label{energyforhol}
\end{equation}
A particularly significant quantity for the solution (\ref{cp2veronese}) satisfying 
the $\field{C}P^2$ model equations (\ref{cp2}) is the topological charge 
\begin{equation}
Q=-\frac{1}{\pi}\int_{S^2} g_{12} d\xi d\bar{\xi}\,,
\label{topologicalcahrge}
\end{equation}
defined on the whole Riemann unit sphere $S^2$. The integral 
(\ref{topologicalcahrge}) exists and is an invariant of the 
surface (\ref{compr81}). It characterizes globally the surface and 
is an integer
\begin{equation}
Q=1\,.
\end{equation}

In the second case, for mixed solutions (\ref{mixedsolcp2}) the 
corresponding orthogonal projector takes the form
\begin{eqnarray}
P_1=\frac{1}{(1+|\xi|^2)^2} \left(\begin{array}{ccc}
(1+|\xi|^4) & -\sqrt{2}\,\xi (|\xi|^2-1) & 2\xi^2 \\
-\sqrt{2}\,\bar{\xi} (|\xi|^2-1) & 4|\xi|^2 & \sqrt{2}\,\xi (|\xi|^2-1) \\
2\bar{\xi}^2 & \sqrt{2}\,\bar{\xi} (|\xi|^2-1)  & 1+|\xi|^4
\end{array}
\right)\,, \label{projectorcp2formix}
\end{eqnarray}
where rank$P_1=2$ and tr$P_1=2$. The surface is determined by 
(\ref{radiusvecmix}) and its induced metric associated with 
the projector (\ref{projectorcp2formix}) is also conformal
\begin{equation}
g_{11}=g_{22}=0, \qquad g_{12}=\frac{2}{(1+|\xi|^2)^2}\,.
\label{inducedmetricformixus}
\end{equation}
The first fundamental form  and the Gaussian curvature are 
\begin{equation}
I=\frac{4}{(1+|\xi|^2)^2}d \xi d \bar{\xi}\,, \qquad \mathcal{K}=1,
\label{firstgaussformixus}
\end{equation}
respectively. 

It is worth mentioning that from the Veronese sequences we can obtain 
associated surfaces with constant Gaussian curvature as stated in 
\cite{Zakrzewski2}. However, the converse statement does not 
apply in general. In this paper and in \cite{ghy} we give examples 
of surfaces with constant Gaussian curvature which are not 
associated with Veronese sequences. Such surfaces correspond to 
the dilation-invariant solutions or mixed soliton solutions of the 
$\field{C}P^{2}$ model.


\section{Comments on surfaces obtained via projector formalism}

Certain geometrical aspects of surfaces 
have been studied recently in \cite{Zakrzewski2}
using the generalized Weierstrass representation 
associated with the $\field{C}P^{N-1}$ sigma models. 
In the context of \cite{Zakrzewski2} a sequence of 
rank-one projectors of the form
\begin{equation}
\mathcal{P}_k:=\mathcal{P}(V_k) = \frac{V_k \otimes V_k^{\dagger}}{V_k^{\dagger} \cdot V_k}\,,
\quad 
\mbox {where}
\,\,\,\,
V_k = {P}_+^k f\,,
\qquad
k \in \field{Z}^+\,,
\label{concproject}
\end{equation}
were used to construct a family of surfaces associated with a given 
solution of the $\field{C}P^{N-1}$ model. 
More specifically, starting with any nonconstant 
holomorphic solution of the $\field{C}P^{N-1}$ model, 
one can successively (say $k$ times) apply the operator $P_+$ 
(given in (\ref{operatorintermsoff})) in order to find a new 
solution $P_k = P_+^k f$, which represents a harmonic map 
$S^2 \rightarrow \field{C}P^{N-1}$. So, for every $k \le N-1$ 
the quantity 
$P_+^k$ constitutes a building element for the construction of a set of 
rank-one projectors 
$\{\mathcal{P}_0, \mathcal{P}_1, \ldots, \mathcal{P}_k\}$ 
as described in (\ref{concproject}). 
This set of projectors determines new conservation laws 
of the form (\ref{conservation}), which 
can be considered new in the sense that only the first one 
is related to the holomorphic 
(or antiholomorphic) solutions and the rest are related 
to the mixed solutions obtained from the nonconstant holomorphic ones by 
applying the operator $P_+$ successively. Consequently, according to this 
procedure it is claimed that one can obtain new surfaces for each projector. 
However, there are some questions concerning this procedure. 
Using the following properties (given in \cite{Zakrzewski})
\begin{eqnarray}
\bar{\partial} ({P}_+^k f) = - {P}_+^{k-1} f \frac{|{P}_+^k f|^2}{|{P}_+^{k-1} f|^2}\,, \nonumber \\
\partial \left(
\frac{{P}_+^{k-1} f}{|{P}_+^{k-1} f|^2} 
\right)
= \frac{{P}_+^k f}{|{P}_+^{k-1} f|^2}\,,
\end{eqnarray}
together with the orthogonality relation (\ref{orthogonalitypr}), it is 
straightforward to compute
\begin{equation}
\partial \mathcal{P}_k = 
\frac{({P}_+^{k+1} f) \otimes ({P}_+^{k} f)^{\dagger}}{|{P}_+^k f|^2} - 
\frac{({P}_+^{k} f) \otimes ({P}_+^{k-1} f)^{\dagger}}{|{P}_+^{k-1} f|^2}
\end{equation}
and 
\begin{equation}
[\partial \mathcal{P}_k, \mathcal{P}_k] = 
\frac{({P}_+^{k+1} f) \otimes ({P}_+^{k} f)^{\dagger}}{|{P}_+^k f|^2} + 
\frac{({P}_+^{k} f) \otimes ({P}_+^{k-1} f)^{\dagger}}{|{P}_+^{k-1} f|^2}
\end{equation}
which can also be written as 
\begin{equation}
[\partial \mathcal{P}_k, \mathcal{P}_k] = 
\partial \mathcal{P}_k + 2 
\frac{({P}_+^{k} f) \otimes ({P}_+^{k-1} f)^{\dagger}}{|{P}_+^{k-1} f|^2}\,.
\label{commutatorwithprank1}
\end{equation}
Similarly, we can write $[\bar{\partial} \mathcal{P}_k, \mathcal{P}_k]$ as 
\begin{equation}
[\bar{\partial} \mathcal{P}_k, \mathcal{P}_k] = 
-\bar{\partial} \mathcal{P}_k - 2 
\frac{({P}_+^{k-1} f) \otimes ({P}_+^{k} f)^{\dagger}}{|{P}_+^{k-1} f|^2}\,.
\label{commutatorwithprank1bar}
\end{equation}
As a consequence of the commutators given in (\ref{commutatorwithprank1}) 
and (\ref{commutatorwithprank1bar}), the Weierstrass data 
(\ref{complex1formimmer}) becomes
\begin{eqnarray}
\!\!dX = -i \left[
(\partial \mathcal{P}_k + 2 
\frac{({P}_+^{k} f) \otimes ({P}_+^{k-1} f)^{\dagger}}{|{P}_+^{k-1} f|^2}) d\xi + 
(\bar{\partial} \mathcal{P}_k + 2 
\frac{({P}_+^{k-1} f) \otimes ({P}_+^{k} f)^{\dagger}}{|{P}_+^{k-1} f|^2}) d\bar{\xi}
\right]\!\!.
\label{weierstrassdataforall}
\end{eqnarray}
It is easily seen that for $k = 0$ (e.g. for the holomorphic solutions, 
or equivalently the antiholomorphic ones) equation 
(\ref{weierstrassdataforall}) reduces to 
\begin{eqnarray}
dX = -i \left[
\partial \mathcal{P}_0 d\xi + 
\bar{\partial} \mathcal{P}_0 d\bar{\xi}
\,\,\right]\,,
\label{weierstrassdataforholantihol}
\end{eqnarray}
since the other terms in (\ref{weierstrassdataforall}) do not appear 
for $k = 0$. Hence, it is concluded that 
$X$ is proportional to the projector $\mathcal{P}_0$. This point has been 
fully discussed both in this paper and in \cite{ghy} and \cite{Hussin}. However, 
for $k \ne 0$ (i.e. for the mixed solutions) the integral of the Weierstrass 
representation (\ref{weierstrassdataforall}) cannot be proportional to 
$\mathcal{P}_k$. 

This point can be further discussed for the example of the mixed solutions 
given in Section 3. This example is also analyzed in 
\cite{Zakrzewski2} by a different approach. 
In \cite{Zakrzewski2} it is stated that 
the mixed solution, obtained from the Veronese vector 
$f=(1, \sqrt{2}\,\xi, \xi^2)$, for the $\field{C}P^{2}$ model 
associated with the projector
\begin{equation}
\mathcal{P}_1=\frac{1}{(1+|\xi|^2)^2} \left(\begin{array}{ccc}
2 |\xi|^2 & \sqrt{2}\,\bar{\xi} (|\xi|^2-1) & -2\bar{\xi}^2 \\
\sqrt{2}\,{\xi} (|\xi|^2-1) & (1-|\xi|^2)^2 & -\sqrt{2}\,\bar{\xi} (|\xi|^2-1) \\
-2{\xi}^2 & -\sqrt{2}\,{\xi} (|\xi|^2-1)  & 2 |\xi|^2
\end{array}
\right)\,, \label{conclusionprojectorcp2formix}
\end{equation}
leads to a radius vector $\vec{{Y}}$ which lies in a $5$-dimensional 
subspace of $\field{R}^{8}$. Moreover, the components of 
the radius vector $\vec{{Y}}$ are given as
\begin{eqnarray}
&&{Y}_1 = \frac{2x(1-x^2-y^2)}{(1+x^2+y^2)^2}\,, \quad
{Y}_2 = \frac{2y(1-x^2-y^2)}{(1+x^2+y^2)^2}\,, \quad
{Y}_3 = \frac{2(x^2-y^2)}{(1+x^2+y^2)^2}\,, \nonumber \\
\!\!\!\!\!\!\!\!\!\!\!\!\!\!\!\!\!\!\!\!
&&{Y}_4 = \frac{4xy}{(1+x^2+y^2)^2}\,, \quad
{Y}_5 = \sqrt{3}\frac{(1-x^2-y^2)^2}{(1+x^2+y^2)^2}\,,
\end{eqnarray}
which satisfy the following surface
\begin{equation}
{Y}_1^2 + {Y}_2^2 + 4 {Y}_3^2 + 4 {Y}_4^2 + 
\frac{1}{\sqrt{3}} {Y}_5 = 1\,.
\end{equation}
However, using the same solution together with the projector 
(\ref{projectorcp2formix}) and the procedure summarized in Section 2, 
we obtain an associated surface with the radius vector $\vec{X}$ 
immersed in a $3$-dimensional subspace of $\field{R}^{8}$. The components 
of the radius vector $\vec{X}$ are given in (\ref{radiusvecmix}) and they 
satisfy equation (\ref{eqsurmix}). Since the two surfaces are obtained 
from the same mixed solutions of the $\field{C}P^{2}$ model, constructed by 
the same procedure from the Veronese vector 
$f=(1, \sqrt{2}\,\xi, \xi^2)$, we expect them to be 
the same geometrical object in accordance with the Bonnet theorem. 
However, it can easily be verified that the two surfaces cannot be 
transformed into each other by rotations and translations.


\section{Dilation-invariant solutions}

The objective of this section is to construct dilation-invariant 
solutions of the $\field{C}P^{2}$ model and then to calculate some 
geometric properties of the surface associated with this model 
by using the Weierstrass formula for immersion in $\field{R}^{8}$. 

Let us discuss the solutions of (\ref{cp2}) which are invariant 
under the scaling symmetries 
\begin{eqnarray}
S = w_i \partial_{w_i}-\bar{w}_i\partial_{\bar{w}_i}\,, \qquad i=1,2\,.
\label{scalsymmcp2}
\end{eqnarray}
For this purpose we determine the invariants of the vector 
fields (\ref{scalsymmcp2}), which imply the algebraic constraints 
$w_i\bar{w}_i = D_i \in \field{R}$, $i=1,2$. 
Without loss of generality we may choose $D_i=1$. Then 
the invariant solution is given by
\begin{equation}
w_i=\frac{F_i(\xi)}{\bar{F}_i(\bar{\xi})}\,, \qquad i=1,2\,,
\label{simplestsolutions}
\end{equation}
where $F_i$ and $\bar{F}_i$ are arbitrary complex-valued functions of 
one complex variable $\xi$ and $\bar{\xi}$, respectively. After 
substituting (\ref{simplestsolutions}) into the $\field{C}P^2$ model 
equations (\ref{cp2}) 
it is immediately seen that the unknown functions 
$F_i$ and $\bar{F}_i$ must satisfy the following differential relation 
\begin{equation}
|F_2|^2 |F_1^{\prime}|^2=|F_1|^2 |F_2^{\prime}|^2\,, 
\label{absconstraint}
\end{equation}
where prime means differentiation with respect to the argument 
({\it i.e.} with respect to either $\xi$ or $\bar{\xi}$). Equation 
(\ref{absconstraint}) implies 
\begin{eqnarray}
F_2^{\prime}({\xi})=
\frac{F_2({\xi})F_1^{\prime}(\xi)}{F_1(\xi)} e^{i\psi}\,,
\label{solabsconstraint}
\end{eqnarray}
which has the following solution
\begin{equation}
F_2({\xi})= c F_1(\xi)^{e^{i\psi}}\,,
\qquad c\in \field{C}\,,
\label{solution}
\end{equation}
where $\psi$ is an arbitrary constant. 
By substituting (\ref{simplestsolutions}) and (\ref{solabsconstraint}) 
into (\ref{cp2}) it is seen that $\psi$ must satisfy
\begin{equation}
\psi=\pm \frac{\pi}{3}+2\pi m\,, \qquad m\in \field{Z}\,.
\label{equationforpsi}
\end{equation}
Thus we obtain a class of scaling invariant solutions of the 
$\field{C}P^2$ model equations (\ref{cp2}) which depend on one 
arbitrary complex-valued function of one variable $\xi$ and its conjugate 
\begin{equation}
w_1=\frac{F_1(\xi)}{\bar{F}_1(\bar{\xi})}\,, \qquad 
w_2=\frac{c}{\bar{c}}\frac{F_1(\xi)^{e^{i\psi}}}
{\bar{F}_1(\bar{\xi})^{e^{-i\psi}}}\,.
\label{generalsolutionsfornonsplitting}
\end{equation}

We now perform a detailed investigation 
of the geometric implications of the induced metric associated with a 
quadratic differential. In \cite{ghy} it was shown that the induced 
metric for the $\field{C}P^2$ model equations subjected to the 
DCs (i.e. $w_i\bar{w}_i = 1$, $i=1,2$) 
is conformal and the Gaussian 
curvature for the associated surfaces vanishes. 
It was also shown that the coordinates of the radius vector $\vec{X}$ 
for the nonsplitting solutions of the $\field{C}P^2$ model equations
are given by  
\begin{eqnarray}
X_1&=&\frac{i}{6\sqrt{3}\,|c|^2}|F|^{-2e^{i\psi}}(\bar{c}^2F-c^2
\bar{F}|F|^{2i\sqrt{3}})\,, \nonumber \\
X_2&=&-\frac{1}{6\sqrt{3}\,|c|^2}|F|^{-2e^{i\psi}}(\bar{c}^2F+c^2
\bar{F}|F|^{2i\sqrt{3}})\,,\nonumber \\
X_3&=&\frac{1}{6}\big((1-i\sqrt{3})\rm{ln}F+(1+i\sqrt{3})\rm{ln}
\bar{F}\big)\,,\nonumber \\
X_4&=&-\frac{1}{6}\big((i+\sqrt{3})\rm{ln}F+(-i+\sqrt{3})\rm{ln}
\bar{F}\big)\,,\nonumber \\
X_5&=&-\frac{F^2+\bar{F}^2}{6\sqrt{3}\, |F|^2}\,,\nonumber \\
X_6&=&\frac{1}{6\sqrt{3}\,|c|^2}|F|^{-2e^{i\psi}}(\bar{c}^2\bar{F}+
c^2{F}|F|^{2i\sqrt{3}})\,,\nonumber \\
X_7&=&\frac{i(F^2-\bar{F}^2)}{6\sqrt{3}\, |F|^2}\,,\nonumber \\
X_8&=&\frac{i}{6\sqrt{3}\,|c|^2}|F|^{-2e^{i\psi}}(\bar{c}^2\bar{F}-
c^2{F}|F|^{2i\sqrt{3}})\,.
\label{rdvecforcompfornonsplitforexamp}
\end{eqnarray}
The corresponding first fundamental form is 
immediately given as
\begin{equation}
I=\frac{2}{3}\, \frac{|F^{\prime}|^2}{|F|^{2}}\,d\xi\,d\bar{\xi}\,.
\label{firstfunforcp2fornonsplitexamp}
\end{equation}
Note that the components of the radius vector $\vec{X}$ in 
(\ref{rdvecforcompfornonsplitforexamp}) satisfy the following relations
\begin{equation}
X_1^2+X_2^2=X_5^2+X_7^2=X_6^2+X_8^2=\frac{1}{27}\,. 
\end{equation}
Eliminating the functions $F$ and $\bar{F}$ in 
(\ref{rdvecforcompfornonsplitforexamp}) we obtain
\begin{eqnarray}
X_1 &=& \frac{i}{6\sqrt{3}\,|c|^2}e^{-(v+\bar{v}) e^{i\psi}}(\bar{c}^2 e^v-c^2 e^{\bar{v}} e^{i \sqrt{3} (v+\bar{v})})\,, \nonumber \\
X_2 &=& -\frac{1}{6\sqrt{3}\,|c|^2}e^{-(v+\bar{v}) e^{i\psi}}(\bar{c}^2 e^v+c^2 e^{\bar{v}} e^{i \sqrt{3} (v+\bar{v})})\,, \nonumber \\
X_5 &=& -\frac{1}{3\sqrt{3}}\cos\left(\frac{3}{2}(\sqrt{3}X_3+X_4)\right)\,, \nonumber \\
X_6 &=& \frac{1}{6\sqrt{3}\,|c|^2}e^{-(v+\bar{v}) e^{i\psi}}(\bar{c}^2 e^{\bar{v}}+c^2 e^{{v}} e^{i \sqrt{3} (v+\bar{v})})\,, \nonumber \\
X_7 &=& -\frac{1}{3\sqrt{3}}\sin\left(\frac{3}{2}(\sqrt{3}X_3+X_4)\right)\,, \nonumber \\
X_8 &=& \frac{i}{6\sqrt{3}\,|c|^2}e^{-(v+\bar{v}) e^{i\psi}}(\bar{c}^2 e^{\bar{v}}-c^2 e^{{v}} e^{i \sqrt{3} (v+\bar{v})})\,,
\end{eqnarray}
where $v=\frac{3}{4} (1+i\sqrt{3}) (X_3+iX_4)$. The surface is parametrized in terms of $X_3$ and $X_4$. Now, the corresponding first fundamental form becomes
\begin{equation}
I=\frac{3}{2}(dX_3^2+dX_4^2)\,.
\end{equation}
Note that this is just the real form of (\ref{firstfunforcp2fornonsplitexamp}) where $\xi^1=X_3$ and $\xi^2=X_4$. 

The induced metric (\ref{firstfunforcp2fornonsplitexamp}) on the $(\xi,\bar{\xi})$ plane can be written as a quadratic differential
\begin{equation}
I=\frac{2}{3}d\left(\ln{F(\xi)}\right)\wedge d\left(\ln{\bar{F}(\bar{\xi})}\right)\,.
\label{Name1}
\end{equation}
Equation (\ref{Name1}) defines a field of line elements on a surface ${\mathcal{F}}$ with singularities at the critical points (i.e. the zeros and poles of the differential (\ref{Name1})). The geodesic trajectories $\xi = \xi(t)$ of this metric are determined locally by the integral
\begin{equation}
\mbox{Re}\left(e^{i\theta}\omega\right) = c\,, \qquad \omega=\int{\frac{F'}{F}}d\xi\,, \qquad \theta \in \field{R}\,,
\end{equation}
where we make use of the definitions and notations given in \cite{Strebel}.
The simplest local trajectory structure of the quadratic differential (\ref{Name1}) can be found by assuming that $F'/F$ has two simple zeros and one simple pole. Then
\begin{equation}
F(\xi)=A\xi^n\left(1+{\mathcal O}(\xi)\right)\,,\quad n\in \field{Z}\,, \qquad A\in \field{C}\,,
\end{equation}
and
\begin{eqnarray}
&&{\frac{F'}{F}}={\frac{n}{\xi}}\left(1+{\mathcal O}(\xi)\right)\mbox{ near a pole at }\xi=0\,, \nonumber \\
&&{\frac{F'}{F}}=C(\xi-a)\left(1+{\mathcal O}(\xi-a)\right)\mbox{ near a simple zero of } F'\,, \quad a \in \field{C}\,.
\end{eqnarray}
Locally, the flat coordinates of the metric $I$ are the real and imaginary parts of the function
\begin{equation}
\omega = \int{\frac{F'}{F}}d\xi=n\ln{\xi}+{\mathcal O}(1)\,,
\label{Name2}
\end{equation}
where ${\mathcal O}(1)$ denotes some analytic function near $\xi=0$. The critical vertical trajectory is defined to be the maximal trajectory of the ODE
\begin{equation}
\mbox{Re}\left({\frac{F'}{F}}d\xi\right) = 0 \qquad\mbox{with }\theta=0\,.
\end{equation}
As a consequence of this equation we obtain the following condition 
\begin{eqnarray}
\frac{d\xi}{dt} = i \frac{\bar{F'}(\xi)}{\bar{F}(\xi)}\,. 
\end{eqnarray}
The monodromy of (\ref{Name2}) is given by $2i\pi n$. Also, the function $q=e^{\omega/n}$ is analytic in a punctured neighborhood of $\xi=0$, since Re$(\omega)=n\ln{|\xi|}$ and $|q|\sim|\xi|$ near $\xi=0$. Thus, $q$ has a removable singularity at $\xi=0$, and hence can be extended to an analytic function
\begin{equation}
q(\xi)=B \xi+{\mathcal O}(\xi^2)\,, \qquad 0 \ne B \in \field{C}\,.
\end{equation}
Let us denote by ${\mathcal D}$ the maximal connected domain foliated by closed trajectories homotopic to a small circle around $\xi=0$. The function $q$ is a single valued conformal map of ${\mathcal D}$ onto the disk of radius $|n|$, since the perimeter of the disk is $2\pi|n|$. Hence, we get one semi-infinite cylinder (homeomorphic to the disk $\{0<|\xi|\leq 1\}$) for each simple pole of $F'/F$. For example, for $F=\xi(\xi-1)$ we have
\begin{equation}
{\frac{F'}{F}}=\frac{2\xi-1}{\xi(\xi-1)} = \frac{1}{\xi} + \frac{1}{\xi-1}\,,\qquad \mbox{Res}_{\infty}
\left({\frac{F'}{F}}\right)=-2.
\end{equation}
Hence we obtain three cylinders, two for the poles at $0$ and $1$, each of perimeter $2\pi$, and one of perimeter $4\pi$ for the pole at $\infty$. If we instead assume that $F'/F$ has two simple zeros and three simple poles, then we get four semi-infinite cylinders glued along the based perimeter with two points of conical singularities. It should be noted that for the $\field{C}P^{1}$ model, 
the surfaces corresponding to the dilation invariant solutions 
are cylinders.\cite{Bracken}

Note that the singular solutions (\ref{generalsolutionsfornonsplitting}) are in fact the meron-like solutions 
(i.e. with logarithmically divergent action at isolated points) of the $\field{C}P^2$ model equations (\ref{cp2}). The meron solutions, obtained from the dilation invariance of the solutions of the 
$\field{C}P^2$ model, are located 
at $\bar{F}(\bar{\xi})=0$. Merons are more durable than instantons in the sense that they can exist in a constant (not necessarily zero) Higgs field. This property is shared by both the $\field{C}P^{N-1}$ and Yang-Mills models \cite{Abbott, Berg, Mihai}.


\section{Conclusions}

In this paper we have shown that if the $\field{C}P^{N-1}$ 
model equations are defined on the sphere $S^2$ and the associated 
action functional of this model is finite, then the specific holomorphic 
function $J$ (i.e. component of the energy-momentum tensor of the model) 
vanishes and consequently the surfaces are conformally 
parametrized. We demonstrate that the holomorphic and the mixed 
Veronese solutions of the $\field{C}P^{2}$ model are associated 
with a sphere and an ellipsoid immersed in $\field{R}^{8}$ and 
$\field{R}^{3}$, respectively. In this context we have shown that the 
Weierstrass formula for immersion of surfaces associated with the mixed 
solutions of the $\field{C}P^{N-1}$ model cannot be proportional to the 
rank-one projectors $\mathcal{P}_k$ for $k \ne 0$, unlike the case $k = 0$. 
The analysis of the geometrical properties of the dilation-invariant 
meron type solutions shows that they represent semi-infinite cylinders.

This research could be expanded in several directions. The meaning of the 
new conservation laws could be investigated in the context of 
surfaces immersed in multi-dimensional spaces. Are they really 
independent? Can they differ from each other by a total divergence 
(as discussed in Section \ref{prilaminary})? Other natural 
directions which could also be addressed involve families of solutions 
obtained recursively and whether they can be related through an 
auto-B\"{a}cklund transformation (note that a B\"{a}cklund parameter 
is not present in formula (\ref{operatorintermsoff}), but it can be introduced 
by a gauge transformation). In addition to these it is also 
important to ask if the symmetry operator 
(\ref{operatorintermsoff}) is expressible in terms of some combination of 
the known infinitesimal generators (as given in \cite{ghy}) of the 
Lie-point symmetry algebra of the $\field{C}P^{N-1}$ model. These and other 
issues will be addressed in our future work.

\section*{Acknowledgments}

The authors would like to thank Marco Bertola (Concordia 
University) and Libor \v{S}nobl (Czech Technical University) 
for helpful and interesting discussions on the topic 
of this paper. The research reported in this paper is supported 
by the NSERC 
of Canada. \.{I}.Y. acknowledges a postdoctoral fellowship 
awarded by the Laboratory of Mathematical Physics of the CRM, 
Universit\'{e} de Montr\'{e}al. This project was completed 
during the A.M.G's visit to the Doppler Institute (Project LC06002 of the Ministry of Education of the Czech Republic) and he would like to thank the 
Institute for their kind invitation.


\end{document}